\documentclass[11pt,a4paper]{amsart}

\usepackage{amsthm,amsmath,graphicx,calc}
\usepackage[lmargin=27mm,rmargin=27mm,tmargin=30mm,bmargin=30mm,footskip=5mm]{geometry}
\usepackage[sort&compress,numbers]{natbib}

\DeclareFontShape{OT1}{cmr}{b}{sc}
{<5><6><7><8><9><10><12><10.95><14.4><17.28><20.74><24.88>cmbcsc10}{}
\DeclareFontShape{OT1}{cmr}{bx}{sc}
{<5><6><7><8><9><10><12><10.95><14.4><17.28><20.74><24.88>cmbcsc10}{}

\DeclareFontShape{OT1}{cmr}{b}{tt}
{<5><6><7><8><9><10><12><10.95><14.4><17.28><20.74><24.88>cmbtt10}{}
\DeclareFontShape{OT1}{cmr}{bx}{tt}
{<5><6><7><8><9><10><12><10.95><14.4><17.28><20.74><24.88>cmbtt10}{}

\newcommand{\aaa}{\textup{(a)}}
\newcommand{\bbb}{\textup{(b)}}
\newcommand{\ccc}{\textup{(c)}}

\newcommand{\half}{\ensuremath{\protect\tfrac{1}{2}}}

\newcommand{\Oh}[1]{\ensuremath{\protect\mathcal{O}(#1)}}

\newcommand{\Figure}[4][htb]{
\begin{figure}[#1]
	\vspace*{1ex}
	\begin{center}#3\end{center}
	\vspace*{-1ex}
	\caption{\figlabel{#2}#4}
\end{figure}
}

\newlength{\marginboxwidth}

\usepackage{pifont}

\newcommand{\thmlabel}[1]{\label{thm:#1}}
\newcommand{\thmref}[1]{Theorem~\ref{thm:#1}}

\newcommand{\lemlabel}[1]{\label{lem:#1}}
\newcommand{\lemref}[1]{Lemma~\ref{lem:#1}}
\newcommand{\twolemref}[2]{Lemmas~\ref{lem:#1} and \ref{lem:#2}}

\newcommand{\figlabel}[1]{\label{fig:#1}}
\newcommand{\figref}[1]{Figure~\ref{fig:#1}}

\newcommand{\seclabel}[1]{\label{sec:#1}}

\theoremstyle{plain}
\newtheorem{theorem}{Theorem}
\newtheorem{lemma}{Lemma}

\renewcommand{\baselinestretch}{1.25}
\begin{document}

\newcommand{\V}[1]{\ensuremath{\protect|#1|}}
\newcommand{\E}[1]{\ensuremath{\protect\|#1\|}}
\newcommand{\CR}[1]{\ensuremath{\protect\textsf{\textup{cr}}(#1)}}
\newcommand{\MCR}[1]{\ensuremath{\protect\textsf{\textup{mcr}}(#1)}}
\newcommand{\RCR}[1]{\ensuremath{\protect\overline{\textsf{\textup{cr}}}(#1)}}
\newcommand{\CCR}[1]{\ensuremath{\protect\textsf{\textup{cr}}^{\star}(#1)}}
\newcommand{\FFF}{\ensuremath{\mathcal{F}}}

\date{\today}
\title[The Minor Crossing Number of Graphs with an Excluded Minor]{The Minor Crossing Number\\ of Graphs with an Excluded Minor}

\author{Drago Bokal}
\address{Department of Combinatorics and Optimization, University of Waterloo, Waterloo, Canada}
\email{dbokal@uwaterloo.ca}

\author{Ga\v{s}per Fijav\v{z}}
\address{Faculty of Computer and Information Science, University of Ljubljana, Ljubljana, Slovenia}
\email{gasper.fijavz@fri.uni-lj.si}

\author{David R. Wood}
\address{Departament de Matem{\'a}tica Aplicada II, Universitat Polit{\`e}cnica de Catalunya, Barcelona, Spain}
\email{david.wood@upc.edu}

\thanks{The research of David Wood is supported by a Marie Curie
Fellowship of the European Community under contract 023865, and by the projects MCYT-FEDER BFM2003-00368 and Gen.\ Cat 2001SGR00224.}


\keywords{graph drawing, crossing number, minor crossing number, graph decomposition, planar decomposition, graph minor}

\subjclass[2000]{05C62 (graph representations), 05C10 (topological graph theory), 05C83 (graph minors)}

\begin{abstract}
The \emph{minor crossing number} of a graph $G$ is the minimum crossing number of a graph that contains $G$ as a minor. It is proved that for every graph $H$ there is a constant $c$, such that every graph $G$ with no $H$-minor has minor crossing number at most $c|V(G)|$.
\end{abstract}


\maketitle

\section{Introduction}
\label{sec:Intro}


The \emph{crossing number} of a graph\footnote{We consider finite, undirected, simple graphs $G$ with vertex set $V(G)$ and edge set $E(G)$. Let $\V{G}:=|V(G)|$ and $\E{G}:=|E(G)|$. Let $\Delta(G)$ be the maximum vertex degree of $G$.} $G$, denoted by \CR{G}, is the minimum number of crossings in a drawing\footnote{A \emph{drawing} of a graph represents each vertex by a distinct point in the plane, and represents each edge by a simple closed curve between its endpoints, such that the only vertices an edge intersects are its own endpoints, and no three edges intersect at a common point (except at a common endpoint). A \emph{crossing} is a point of intersection between two edges (other than a common endpoint). A graph is \emph{planar} if it has a crossing-free drawing.} of $G$ in the plane; see \citep{PachToth-Geom00, PachToth-JCTB00, Vrto, Szekely-DM04, Szekely-SOFSEM05, EG-AMM73, RS-CrossingNumberSurvey} for surveys. The crossing number is an important measure of the non-planarity of a graph \citep{Szekely-DM04}, with applications in discrete and computational geometry \citep{Szekely-CPC97, PachSharir-CPC98} and VLSI circuit design \citep{BL84, Leighton84, Leighton83}. In information visualisation, one of the most important measures of the quality of a graph drawing is the number of crossings \citep{Purchase-JVLC98, Purchase97, PCJ97}. 

We now outline various aspects of the crossing number that have been studied. First note that computing the crossing number is $\mathcal{NP}$-hard \citep{GJ-SJDM83}, and remains so for simple cubic graphs \citep{PSS05, Hliney-JCTB06}. Moreover, the exact or even asymptotic crossing number is not known for specific graph families, such as complete graphs \citep{RT-AMM97}, complete bipartite graphs \citep{RS-JGT96, Nahas-EJC03, RT-AMM97}, and cartesian products \citep{BP, BT, AR-JCTB04, RT-DCG95, GS-JGT04}. Given that the crossing number seems so difficult, it is natural to focus on asymptotic bounds rather than exact values. The `crossing lemma', conjectured by \citet{EG-AMM73} and 
first proved by \citet{Leighton83} and \citet{Ajtai82}, gives such a lower bound. It states that for some constant $c$, $\CR{G}\geq c\E{G}^3/\V{G}^2$ for every graph $G$ with $\E{G}\geq4\V{G}$. 
See \citep{PRTT-SoCG04, Montaron-JGT05} for recent improvements. Other general lower bound techniques that arose out of the work of \citet{Leighton83, Leighton84} include the bisection/cutwidth method \citep{PSS-Algo96, SSSV-Algo96, DV-JGAA03, SS-CPC94} and the embedding method \citep{SS-CPC94, SSSV-AM96}. Upper bounds on the crossing number of general families of graphs have been less studied. One example, by \citet{PachToth-GD05}, says that graphs $G$ of bounded genus and bounded degree have \Oh{\V{G}} crossing number. See \citep{BPT-IJFCS06, DjidVrto-ICALP06} for extensions. The present paper also focuses on crossing number upper bounds.

Graph minors\footnote{Let $vw$ be an edge of a graph $G$.
Let $G'$ be the graph obtained by identifying the vertices $v$ and $w$, deleting
loops, and replacing parallel edges by a single edge. Then $G'$ is obtained from
$G$ by \emph{contracting} $vw$. A graph $H$ is a \emph{minor} of a graph $G$ if
$H$ can be obtained from a subgraph of $G$ by contracting edges. A family of
graphs \FFF\ is \emph{minor-closed} if $G\in\FFF$ implies that every minor of
$G$ is in \FFF. \FFF\ is \emph{proper} if it is not the family of all graphs. A
deep theorem of \citet{RS-GraphMinorsXX-JCTB04} states that every proper
minor-closed family can be characterised by a finite family of excluded minors.
Every proper minor-closed family is a subset of the $H$-minor-free graphs for
some graph $H$. We thus focus on minor-closed families with one excluded minor.} are a widely used structural tool in graph theory. So it is inviting to explore the relationship between minors and the crossing number. One impediment
is that the crossing number is not minor-monotone; that is, there are graphs $G$ and $H$ with $H$ a minor of $G$, for which $\CR{H}>\CR{G}$. Nevertheless, following an initial paper by \citet{RS-GST93}, there have been a number of recent papers on the relationship between crossing number and graph minors \citep{GS-JGT01, BFM-SJDM06, PSS05, GRS-EJC04, Negami-JGT01, BCSV-ENDM, Hliney-JCTB03, Hliney-JCTB06, WoodTelle}. For example, \citet{WoodTelle} proved the following upper bound (generalising the above-mentioned results in \citep{PachToth-GD05, BPT-IJFCS06, DjidVrto-ICALP06} for graphs of bounded genus).

\begin{theorem}[\citep{WoodTelle}]
\thmlabel{CrossingNumber}
For every graph $H$ there is a constant $c=c(H)$, such that every $H$-minor-free graph $G$ has crossing number $\CR{G}\leq c\,\Delta(G)^2\V{G}$.
\end{theorem}

\subsection{Minor Crossing Number}

\citet{BFM-SJDM06} defined the \emph{minor crossing number} of a graph $G$, denoted by \MCR{G}, to be the minimum crossing number of a graph that contains $G$ as a minor. The main motivation for this definition is that for every constant $c$, the family of graphs $G$ for which $\MCR{G}\leq c$ is closed under taking minors. Moreover, the minor crossing number corresponds to a natural style of graph drawing, in which each vertex is drawn as a tree. \citet{BCSV-ENDM} proved a number of lower bounds on the minor crossing number that parallel the lower bound techniques of Leighton. The main result of this paper is to prove the following upper bound, which is an analogue of \thmref{CrossingNumber} for the minor crossing number (without the dependence on the maximum degree). 

\begin{theorem}
\thmlabel{MinorCrossingNumber}
For every graph $H$ there is a constant $c=c(H)$, such that every $H$-minor-free graph $G$ has minor crossing number $\MCR{G}\leq c\,\V{G}$.
\end{theorem}

The restriction to graphs with an excluded minor in \thmref{MinorCrossingNumber} is unavoidable in the sense that $\MCR{K_n}\in\Theta(n^2)$. The linear dependence in \thmref{MinorCrossingNumber} is best possible since  $\MCR{K_{3,n}}\in\Theta(n)$. Both these bounds were established by  \citet{BFM-SJDM06}. An elegant feature of \thmref{MinorCrossingNumber} and the minor crossing number is that there is no dependence on the maximum degree,  unlike in  \thmref{CrossingNumber}, where some dependence on the maximum degree is unavoidable. In particular, the complete bipartite graph $K_{3,n}$ has no $K_5$-minor and has $\Theta(n^2)$ crossing number \citep{RS-JGT96, Nahas-EJC03}.

\section{Planar Decompositions}
\seclabel{Decompositions}

It is widely acknowledged that the theory of crossing numbers needs new ideas. Some tools that have been recently developed include `meshes'  \citep{RT-DCG95}, `arrangements' \citep{AR-JCTB04}, `tile drawings' \citep{PR-JGT03, PR-JGT03, BCCG, PR-AJC04}, and the `zip product' \citep{BCCG, BP, BT}. A feature of the proof of \thmref{CrossingNumber} by \citet{WoodTelle} is the use of `planar decompositions' as a new tool for studying the crossing number. Planar decompositions are also the key component in the proof of \thmref{MinorCrossingNumber} in this paper. 

Let $G$ and $D$ be graphs, such that each vertex of $D$ is a set of vertices of $G$ (called a \emph{bag}). Note that we allow distinct vertices of $D$ to be the same set of vertices in $G$; that is, $V(D)$ is a multiset. For each vertex $v$ of $G$, let $D(v)$ be the subgraph of $D$ induced by the bags that contain $v$. Then $D$ is a \emph{decomposition} of $G$ if:
\begin{itemize}
\item $D(v)$ is connected and nonempty for each vertex $v$ of $G$, and 
\item $D(v)$ and $D(w)$ touch\footnote{Let $A$ and $B$ be subgraphs of a graph $G$. Then $A$ and $B$ \emph{intersect} if $V(A)\cap V(B)\ne\emptyset$, and $A$ and $B$ \emph{touch} if they intersect or $v\in V(A)$ and $w\in V(B)$ for some edge $vw$ of $G$. } for each edge $vw$ of $G$.
\end{itemize}

Decompositions, when $D$ is a tree, were first studied in detail by \citet{RS-GraphMinors-II}. \citet{DiestelKuhn-DAM05}\footnote{A decomposition was called a \emph{connected decomposition} by \citet{DiestelKuhn-DAM05}.}  first generalised the definition for arbitrary graphs $D$.

We measure the `complexity' of a graph decomposition $D$ by the following parameters. The \emph{width} of $D$ is the maximum cardinality of a bag. The \emph{order} of $D$ is the number of bags. The \emph{degree} of $D$ is the maximum degree of the graph $D$. The decomposition $D$ is \emph{planar} if the graph $D$ is planar. 


\citet{DiestelKuhn-DAM05} observed that decompositions generalise minors in the following sense. 

\begin{lemma}[\citep{DiestelKuhn-DAM05}]
\lemlabel{Minor2Decomp}
A graph $G$ is a minor of a graph $D$ if and only if a graph isomorphic to $D$ is a decomposition of $G$ with width $1$.
\end{lemma}


\citet{WoodTelle} describe a number of tools for manipulating decompositions, such as the following lemma for composing two decompositions.

\begin{lemma}[\citep{WoodTelle}]
\lemlabel{Composition}
Suppose that $D$ is a decomposition of a graph $G$ with width $k$, and that $J$ is a decomposition of $D$ with width $\ell$. Then $G$ has a decomposition isomorphic to $J$ with width $k\ell$. 
\end{lemma}

\lemref{Composition} has the following special case.

\begin{lemma}
\lemlabel{SpecialComposition}
If a graph $G_1$ is a minor of a graph $G_2$, and $J$ is a decomposition of $G_2$ with width $\ell$, then a graph isomorphic to $J$ is a decomposition of $G_1$ with width $\ell$. 
\end{lemma}

\begin{proof}
By \lemref{Minor2Decomp}, a graph isomorphic to $G_2$ is a decomposition of $G_1$ with width $1$. By \lemref{Composition} with $k=1$ and $D=G_2$, a graph isomorphic to $J$ is a decomposition of $G_1$ with width $\ell$. 
\end{proof}

The next tool by \citet{WoodTelle} reduces the order of a decomposition at the expense of increasing the width.

\begin{lemma}[\citep{WoodTelle}]
\lemlabel{ReduceOrder}
Suppose that a graph $G$ has a planar decomposition $D$ of width $k$ and order at most $c\V{G}$ for some $c\geq1$. Then $G$ has a planar decomposition of width $c'k$ and order $\V{G}$, for some $c'$ depending only on $c$. 
\end{lemma}


Converse to \lemref{ReduceOrder}, we now show that the width and degree of a decomposition can be reduced at the expense of increasing the order.

\begin{lemma}
\lemlabel{Simplify}
If a graph $G$ has a planar decomposition $D$ of width $k$, then $G$ has:\\
\aaa\ a planar decomposition $D_1$ of width $k$, order $\V{D_1}<6\V{D}$, and degree $\Delta(D_1)\leq 3$, \\
\bbb\ a planar decomposition $D_2$ of width $2$, order $\V{D_2}<3k(k+1)\V{D}$, and degree $\Delta(D_2)\leq 4$,\\
\ccc\ a planar decomposition $D_3$ of width $2$, order $\V{D_3}<6k^2\V{D}$, and degree $\Delta(D_3)\leq 3$.
\end{lemma}

\begin{proof} 
Fix an embedding of $D$ in the plane. By adding edges we can assume that $D$ has minimum degree at least $3$. 

First we prove (a). Let $D_1$ be the graph with two vertices $X_e$ and $Y_e$ for every edge $e=XY\in E(D)$, where each bag $X_e$ is a copy of $X$. We say that $X_e$ \emph{belongs} to $X$. Add the edge $X_eY_e$ to $D_1$ for each edge $e=XY\in E(D)$. Add the edge $X_eX_f$ to $D_1$ whenever the edges $e$ and $f$ are consecutive in the cyclic order of edges incident to a bag $X$ in $D$. 

As illustrated in \figref{Simplify}(b), each bag $X$ is thus replaced by a cycle in $D_1$, each vertex of which has one more incident edge in $D_1$. Thus $D_1$ is a planar graph with maximum degree $3$ and order $\V{D_1}=2\E{D}$ (after adding edges to $D$). Since $D$ is planar, $\E{D}\leq 3\V{D}-6$ and $\V{D_1}\leq 6\V{D}-12$. Since the set of bags of $D_1$ that belong to a specific bag of $D$ induce a connected (cycle) subgraph of $D_1$, and $D(v)$ is a connected subgraph of $D$ for each vertex $v$ of $G$, $D_1(v)$ is a connected subgraph of $D_1$. 

We now prove that $D_1(v)$ and $D_1(w)$ touch for each edge $vw$ of $G$. If $v$ and $w$ are in a common bag $X$ of $D$, then $v$ and $w$ are in every bag $X_e$ of $D_1$. Otherwise, $v\in X$ and $w\in Y$ for some edge $e=XY$ of $D$, in which case $v\in X_e$, $w\in Y_e$, and $X_eY_e$ is an edge of $D_1$. Thus $D_1(v)$ and $D_1(w)$ touch. Therefore $D_1$ is a planar decomposition of $G$. This completes the proof of (a). 

\Figure[h]{Simplify}{\includegraphics{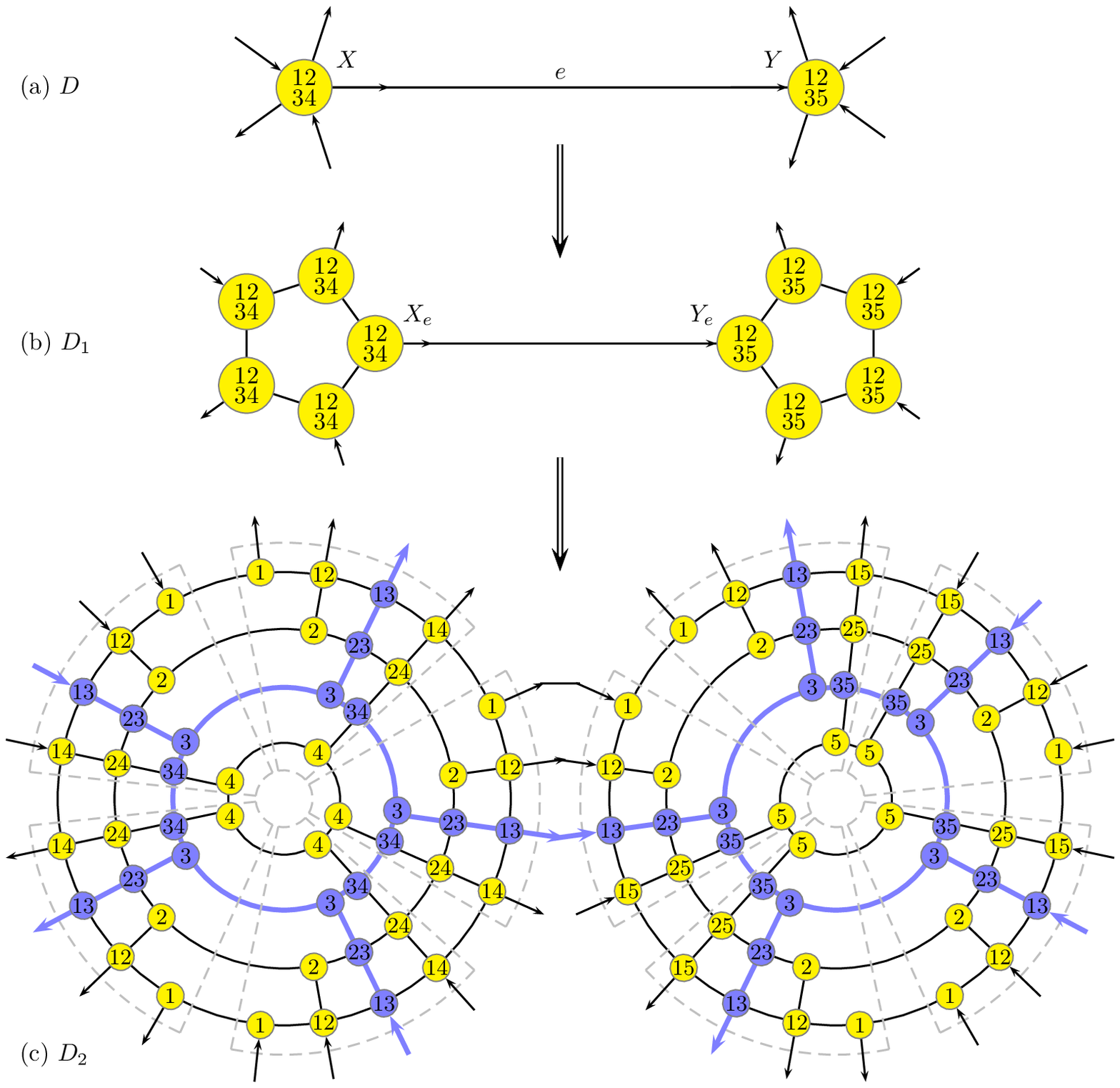}}{(a) The planar decomposition $D$. (b) The planar decomposition $D_1$ obtained from $D$ by replacing each bag of degree $d$ by $d$ bags of degree $3$. (c) The planar decomposition $D_2$ obtained from $D_1$ by replacing each bag of width $k$ by $\binom{k+1}{2}$ bags of width $2$. The subgraph $D_2(3)$ is highlighted.}

Now we prove (b). Fix an arbitrary linear order $\preceq$ on $V(G)$, and arbitrarily orient the edges of $D$. For each arc $e=\overrightarrow{XY}$ of $D$, orient the edge $X_eY_e$ of $D_1$ from $X_e$ to $Y_e$. Informally speaking, we now construct a planar decomposition $D_2$ from $D_1$ by replacing each bag $X_e$ of $D_1$ by a set of $\binom{k+1}{2}$ bags, each of width $1$ or $2$, that form a wedge pattern, as illustrated in \figref{Simplify}(c). Depending on whether $e$ is incoming or outgoing at $X$, the wedge is reflected appropriately to ensure the planarity of $D_2$.

Now we define $D_2$ formally. Consider a bag $X=\{v_1,v_2,\dots,v_k\}$ of $D$, where $v_1\prec v_2\prec \dots\prec v_k$. For each pair of vertices $v_i,v_j$ in $X$, and for each edge $e$ incident to $X$, add a bag labelled $\{v_i,v_j\}_{X_e}$ to $D_2$, where $\{v_i,v_j\}_{X_e}$ is a copy of $\{v_i,v_j\}$. (The bag $\{v_i,v_i\}_{X_e}$ is a singleton $\{v_i\}$.)\ We say that $\{v_i,v_j\}_{X_e}$ \emph{belongs} to $X_e$ and to $X$. Thus there are $\binom{k+1}{2}$ bags that belong to each bag of $D_1$. Hence $\V{D_2}\leq\binom{k+1}{2}\V{D_1}<3k(k+1)\V{D}$. Add an edge in $D_2$ between the bags $\{v_i,v_j\}_{X_e}$ and $\{v_i,v_{j+1}\}_{X_e}$ for $1\leq i\leq k$ and $1\leq j\leq k-1$. As illustrated in \figref{Simplify}(c), the subgraph of $D_2$ induced by the bags that belong to each bag $X_e$ of $D_1$ form a planar grid-like graph.

Consider two edges $e=XY$ and $f=XZ$ of $D$ that are consecutive in the cyclic order of edges incident to a bag $X$ of $D$ (defined by the planar embedding). Since $D$ has minimum degree $3$, without loss of generality, $XZ$ is clockwise from $XY$. We now add edges to $D_2$ between certain bags that belong to $X_e$ and $X_f$ depending on the orientations of the edges $XY$ and $XZ$. For $1\leq i\leq k$, let $P_i$ be the bag $\{v_i,v_k\}_{X_e}$ if $\overrightarrow{XY}$ and $\{v_i,v_i\}_{X_e}$ if $\overrightarrow{YX}$, and let $Q_i$ be the bag $\{v_i,v_i\}_{X_f}$ if $\overrightarrow{XZ}$ and $\{v_i,v_k\}_{X_f}$ if $\overrightarrow{ZX}$. Add an edge between $P_i$ and $Q_i$ for $1\leq i\leq k$. As illustrated in \figref{Simplify}(c), the subgraph of $D_2$ induced by the bags that belong to each bag $X$ of $D$ is planar. 

Now consider an edge $e=\overrightarrow{XY}$ of $D$, where $X=\{v_1,v_2,\dots,v_k\}$ with $v_1\prec v_2\prec \dots\prec v_k$, and $Y=\{w_1,w_2,\dots,w_k\}$ with $w_1\prec w_2\prec \dots\prec w_k$. Whenever $v_i=w_j$, add an edge between $\{v_1,v_i\}_{X_e}$ and $\{w_1,w_j\}_{Y_e}$ to $D_2$. This completes the construction of $D_2$. Observe that the bags $\{v_1,v_1\}_{X_e},\{v_1,v_2\}_{X_e},\dots,\{v_1,v_k\}_{X_e}$ are ordered clockwise on the outer face of the subgraph of $D_2$ induced by the bags belonging to $X$. Similarly, the bags $\{w_1,w_1\}_{Y_e},\{w_1,w_2\}_{Y_e},\dots,\{w_1,w_k\}_{Y_e}$ are ordered anticlockwise on the outer face of the subgraph of $D_2$ induced by the bags belonging to $Y$. Thus these edges do not introduce any crossings in $D_2$, as illustrated in \figref{Simplify}(c).

We now prove that each subgraph $D_2(v)$ is a nonempty connected subgraph of $D_2$ for each vertex $v$ of $G$. Say $v$ is in a bag $X=\{v_1,v_2,\dots,v_k\}$ of $D$, with $v_1\prec v_2\prec \dots\prec v_k$ and $v=v_i$. Observe that the set of bags $\{\{v_i,v_j\}_{X_e}:v_j\in X\in e\in E(D),i\leq j\}$ form a cycle in $D_2$ (drawn as a circle in \figref{Simplify}(c)), and for each edge $e$ incident to $X$, the bags $\{\{v_i,v_j\}_{X_e}:v_j\in X\in e\in E(D),j\leq i\}$ form a path between $\{v_1,v_i\}_{X_e}$ and $\{v_i,v_i\}_{X_e}$, where it attaches to this cycle. Thus the set of bags in $D_2$ that belong to $X$ and contain $v$ form a connected subgraph of $D_2$. For each edge $e=XY$ of $D$ with $v\in X\cap Y$, there is an edge in $D_2$ (between some bags $\{v_1,v_i\}_{X_e}$ and $\{w_1,w_j\}_{Y_e}$) that connects the the set of bags that belong to $X$ and contain $v$ with the set of bags that belong to $Y$ and contain $v$. Thus $D_2(v)$ is connected since $D(v)$ is connected.

We now prove that $D_2(v)$ and $D_2(w)$ touch for each edge $vw$ of $G$. If $v$ and $w$ are in a common bag $X$ of $D$, then $v$ and $w$ are in every bag $\{v,w\}_{X_e}$ of $D_1$. Otherwise, $v\in X$ and $w\in Y$ for some edge $e=XY$ of $D$, in which case $v$ and $w$ are in adjacent bags $\{v_1,v\}_{X_e}$ and $\{w_1,w\}_{Y_e}$, for appropriate vertices $v_1$ and $w_1$. Thus $D_2(v)$ and $D_2(w)$ touch. Therefore $D_2$ is a decomposition of $G$. Observe that $\Delta(D_2)\leq 4$. This completes the proof of (b). 

Now we prove (c). Construct a planar decomposition $D_3$ from $D_2$ by the following operation applied to each bag $X$ of $D_2$ with degree $4$. Say the neighbours of $X$ are $Y_1,Y_2,Y_3,Y_4$ in clockwise order in the embedding of $D_2$. Replace $X$ by two bags $X_1$ and $X_2$, both copies of $X$, where $X_1$ is adjacent to $X_2,Y_1,Y_2$, and $X_2$ is adjacent to $X_1,Y_3,Y_4$. Clearly $D_3$ is a planar decomposition of $G$ with maximum degree $3$. For each bag $X_e$ of $D_1$, there are $k$ bags of degree $3$ in $D_2$ and $\half k(k-1)$ bags of degree $4$ that belong to $X_e$. Since each bag of degree $4$ in $D_2$ is replaced by two bags in $D_3$, there are $k+2(\half k(k-1))=k^2$ bags in $D_3$ that belong to $X_e$. Thus $\V{D_3}\leq 2k^2\E{D}<6k^2\V{D}$.  This completes the proof of (c). 
\end{proof}


Note that the upper bound of $\V{D_1}\leq 6\V{D}$ in \lemref{Simplify}(a) can be improved to $\V{D_1}\leq 4\V{D}$ by replacing each bag of degree $d$ by $d-2$ bags of degree $3$, as illustrated in \figref{MakeSmallCubic}. We omit the details.



\Figure{MakeSmallCubic}{\includegraphics{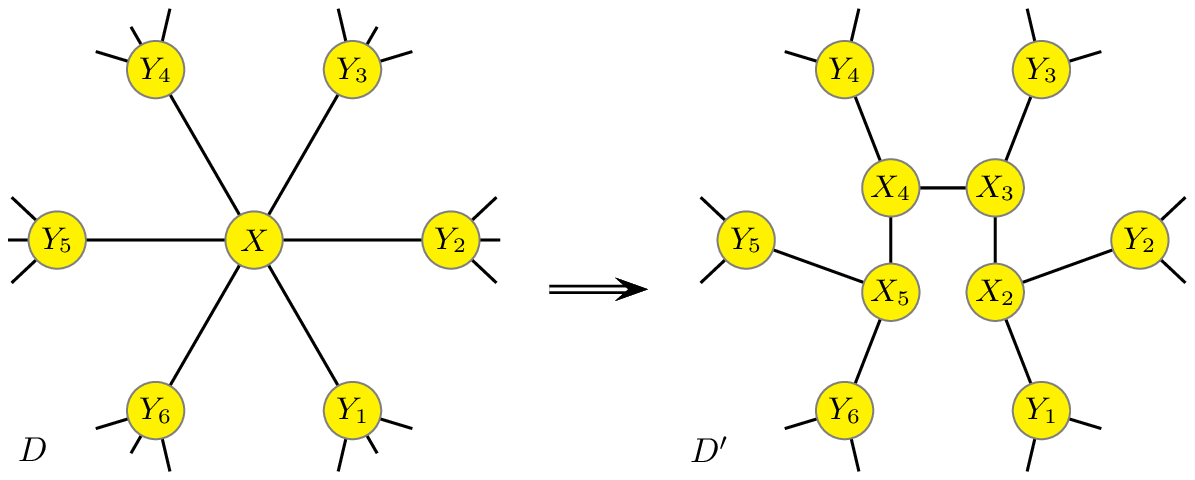}}{Replacing a bag of degree $6$ by four bags of degree $3$.}








\section{Planar Decompositions and Crossing Number}

In this section we review some of the results by \citet{WoodTelle} that link planar decompositions and crossing number.

\begin{lemma}[\citep{WoodTelle}]
\lemlabel{DecompToCrossing}
If $D$ is a planar decomposition of a graph $G$ with width $k$, then $G$ has crossing number
\begin{equation*}
\CR{G}
\;\leq\;k(k+1)\,\Delta(G)^2\,\V{D}\enspace.
\end{equation*}
\end{lemma}

\begin{lemma}[\citep{WoodTelle}]
\lemlabel{DecompMinorFree}
For every graph $H$ there is an integer $k=k(H)$, such that every $H$-minor-free graph $G$ has a planar decomposition of width $k$ and order $\V{G}$.
\end{lemma}

Observe that \twolemref{DecompToCrossing}{DecompMinorFree} imply \thmref{CrossingNumber}. The next lemma is converse to \lemref{DecompToCrossing}.

\begin{lemma}[\citep{WoodTelle}]
\lemlabel{CrossingToDecomp}
Every graph $G$ has a planar decomposition of width $2$ and order $\V{G}+\CR{G}$.
\end{lemma}

We have the following characterisation of graphs with linear crossing number. 

\begin{theorem}[\citep{WoodTelle}]
\thmlabel{CrossingNumberCharacterisation}
The following are equivalent for a graph $G$ of bounded degree:
\begin{enumerate}
\item $\CR{G}\leq c_1\V{G}$ for some constant $c_1$,
\item $G$ has a planar decomposition with width $c_2$ and order $\V{G}$ for some constant $c_2$,
\item $G$ has a planar decomposition with width $2$ and order $c_3\V{G}$ for some constant $c_3$.
\end{enumerate}
\end{theorem}

\begin{proof}
\lemref{CrossingToDecomp} implies that (1) $\Rightarrow$ (3). 
\lemref{ReduceOrder} implies that (3) $\Rightarrow$ (2).
\lemref{DecompToCrossing} implies that (2) $\Rightarrow$ (1).
\end{proof}

Note that \lemref{Simplify}(c) provides a more direct proof that (2) $\Rightarrow$ (3) in \thmref{CrossingNumberCharacterisation} (without the dependence on degree).

\section{Planar Decompositions and Minor Crossing Number}

\lemref{DecompToCrossing} can be extended to give the following upper bound on the minor crossing number. Basically we replace the dependence on $\Delta(G)$ in \lemref{DecompToCrossing} by $\Delta(D)$. 

\begin{lemma}
\lemlabel{DecompToMinorCrossing}
If $D$ is a planar decomposition of a graph $G$ with width $k$, then $G$ has minor crossing number
\begin{equation*}
\MCR{G}\;<\;k^3(k+1)(\Delta(D)+1)^2\,\V{D}\enspace.
\end{equation*}
\end{lemma}

\begin{proof}
Let $G'$ be the graph with one vertex for each occurrence of a vertex of $G$ in a bag of $D$. Consider a vertex $x$ of $G'$ in bag $X$ and a distinct vertex $y$ of $G'$ in bag $Y$. Connect $x$ and $y$ by an edge in $G'$ if and only if $X=Y$ or $XY$ is an edge of $D$. ($G'$ is a subgraph of the lexicographic product $D[K_k]$.)\ For each vertex $v$ of $G$, the copies of $v$ form a connected subgraph of $G'$, since $D(v)$ is a connected subgraph of $D$. Since $D(v)$ and $D(w)$ touch for each edge $vw$ of $G$, some copy of $v$ is adjacent to some copy of $w$. Thus $G$ is a minor of $G'$, and $\MCR{G}\leq\CR{G'}$. Moreover, $D$ defines a planar decomposition of $G'$ with width $k$. By \lemref{DecompToCrossing} applied to $G'$, 
\begin{equation*}
\MCR{G}
\;\leq\;
\CR{G'}
\;\leq\;k(k+1)\,\Delta(G')^2\,\V{D}\enspace.
\end{equation*}
A neighbour of a vertex $x$ of $G'$ is in the same bag as $x$ or is in a neighbouring bag. Thus $\Delta(G')\leq(\Delta(D)+1)k-1$. Thus
\begin{equation*}
\MCR{G}
\;<\;k(k+1)\,((\Delta(D)+1)k)^2\,\V{D}
\;=\;k^3(k+1)\,(\Delta(D)+1)^2\,\V{D}
\enspace.
\end{equation*}
\end{proof}


\twolemref{DecompToMinorCrossing}{Simplify}(a) imply that if $D$ is a planar decomposition of a graph $G$ with width $k$, then $G$ has minor crossing number in \Oh{k^4\V{D}}. This bound can be improved by further transforming the decomposition into a decomposition with width $2$. In particular, \twolemref{DecompToMinorCrossing}{Simplify}(b) imply:

\begin{lemma}
\lemlabel{Decomp2MCR}
If $D$ is a planar decomposition of a graph $G$ with width $k$, then $G$ has minor crossing number
\begin{equation*}
\MCR{G}\;<\;2^3(2+1)(4+1)^2\,3k(k+1)\V{D}\;=\;1800\,k(k+1)\V{D}
\enspace.
\end{equation*}
\end{lemma}

Observe that \twolemref{DecompMinorFree}{Decomp2MCR} imply \thmref{MinorCrossingNumber}. We now set out to prove a converse result to \thmref{MinorCrossingNumber}. 

\begin{lemma}
\lemlabel{SmallRealiser}
For every graph $G$, there is a graph $G'$ containing $G$ as a minor, such that $\MCR{G}=\CR{G'}$ and $\V{G'}\leq\V{G}+\MCR{G}$. 
\end{lemma}

\begin{proof}
By definition, there is a graph $G'$ containing $G$ as a minor, such that $\MCR{G}=\CR{G'}$. Choose such a graph $G'$ with the minimum number of vertices. There is a set $\{T_v:v\in V(G)\}$ of disjoint subtrees in $G'$, such that for every edge $vw$ of $G$, some vertex of $T_v$ is adjacent to some vertex of $T_w$. Every vertex of $G'$ is in some $T_v$, as otherwise we could delete the vertex from $G'$. Hence
\begin{equation*}
\V{G'}
\;=\;\sum_{v\in V(G)}\V{T_v}
\;=\;\V{G}+\sum_{v\in V(G)}(\V{T_v}-1)
\;=\;\V{G}+\sum_{v\in V(G)}\E{T_v}\enspace.
\end{equation*}
We can assume that every edge of every subtree $T_v$ is in some crossing, as otherwise we could contract the edge. Thus
$\V{G'}
\,\leq\,\V{G}+\CR{G'}
\,=\,\V{G}+\MCR{G}$.
\end{proof}

The next lemma is an analogue of \lemref{CrossingToDecomp}.

\begin{lemma}
\lemlabel{MCR2Decomp}
Every graph $G$ has a planar decomposition with width $2$ and order $\V{G}+2\,\MCR{G}$.
\end{lemma}

\begin{proof}
By \lemref{SmallRealiser}, there is some graph $G'$ containing $G$ as a minor, such that $\CR{G'}=\MCR{G}$ and $\V{G'}\leq\V{G}+\MCR{G}$. By \lemref{CrossingToDecomp}, $G'$ has a planar decomposition of width $2$ and order $\V{G'}+\CR{G'}=\V{G'}+\MCR{G}\leq \V{G}+2\,\MCR{G}$. By \lemref{SpecialComposition}, $G$ has a planar decomposition with the same properties.
\end{proof}


We have the following characterisation of graphs with linear minor crossing number, which is analogous to \thmref{CrossingNumberCharacterisation} for crossing number (without the dependence on degree).

\begin{theorem}
The following are equivalent for a graph $G$:
\begin{enumerate}
\item $\MCR{G}\leq c_1\V{G}$ for some constant $c_1$,
\item $G$ has a planar decomposition with width $c_2$ and order $\V{G}$ for some constant $c_2$,
\item $G$ has a planar decomposition with width $2$ and order $c_3\V{G}$ for some constant $c_3$.
\end{enumerate}
\end{theorem}

\begin{proof}
\lemref{MCR2Decomp} implies (1) $\Rightarrow$ (3). 
\lemref{ReduceOrder} implies that (3) $\Rightarrow$ (2).
\lemref{Decomp2MCR} implies that (2)  $\Rightarrow$ (1). 
\end{proof}


\def\soft#1{\leavevmode\setbox0=\hbox{h}\dimen7=\ht0\advance \dimen7
  by-1ex\relax\if t#1\relax\rlap{\raise.6\dimen7
  \hbox{\kern.3ex\char'47}}#1\relax\else\if T#1\relax
  \rlap{\raise.5\dimen7\hbox{\kern1.3ex\char'47}}#1\relax \else\if
  d#1\relax\rlap{\raise.5\dimen7\hbox{\kern.9ex \char'47}}#1\relax\else\if
  D#1\relax\rlap{\raise.5\dimen7 \hbox{\kern1.4ex\char'47}}#1\relax\else\if
  l#1\relax \rlap{\raise.5\dimen7\hbox{\kern.4ex\char'47}}#1\relax \else\if
  L#1\relax\rlap{\raise.5\dimen7\hbox{\kern.7ex
  \char'47}}#1\relax\else\message{accent \string\soft \space #1 not
  defined!}#1\relax\fi\fi\fi\fi\fi\fi} \def\cprime{$'$}
  \def\Dbar{\leavevmode\lower.6ex\hbox to 0pt{\hskip-.23ex \accent"16\hss}D}

\end{document}